\documentclass[12pt]{article}

\usepackage{amsthm}
\usepackage[dvips]{graphicx}
\usepackage{epsfig}
\usepackage{latexsym}
\usepackage{amsmath}
\usepackage{amssymb}
\usepackage{amsfonts}
\usepackage{amscd}

\def\no{\noindent}
\theoremstyle{definition}

\newtheorem{dfn}{Definition}[section]

\newtheorem{ex}[dfn]{Example}

\newtheorem{rem}[dfn]{Remark}

\theoremstyle{plain}

\newtheorem{thm}[dfn]{Theorem}
\newtheorem{mthm}[dfn]{Main Theorem}
\newtheorem{lem}[dfn]{Lemma}
\newtheorem{transferthm}[dfn]{Transfer Theorem}

\newtheorem{prop}[dfn]{Proposition}

\newtheorem{add}[dfn]{Addendum}

\def\proof{\par\medskip\noindent{\it Proof: }}

\def\lra{\longrightarrow}

\def\ra{\rightarrow}

\def\acts{\curvearrowright}

\def\embed{\hookrightarrow}

\def\R{{\mathbb R}}

\def\Z{{\mathbb Z}}

\def\eps{\epsilon}
\def\al{\alpha}

\def\ga{\gamma}

\def\de{\delta}
\def\De{\Delta}
\def\Si{\Sigma}
\def\si{\sigma}

\def\3{\ss}
\def\8{\infty}

\def\geo{\partial_{\infty}}
\def\half{\frac{1}{2}}

\def\tangle{\angle_{Tits}}

\def\ol{\overline}

\def\tits{\partial_{Tits}}

\def\ola{\overrightarrow}

\begin{document}

\title{Polygons in buildings and their refined side lengths}
\author{Misha Kapovich, Bernhard Leeb and John Millson}
\date{November 6, 2005}

\maketitle

\begin{abstract}
\no
As in a symmetric space of noncompact type,
one can associate to an oriented geodesic segment in a Euclidean building 
a vector valued length in the Euclidean Weyl chamber $\De_{euc}$.
In addition to the metric length
it contains information on the direction of the segment. 
We study in this paper restrictions on the $\De_{euc}$-valued 
side lengths of polygons in Euclidean buildings. 
The main result is that 
for thick Euclidean buildings $X$ 
the set ${\cal P}_n(X)$ 
of possible $\De_{euc}$-valued 
side lengths of oriented $n$-gons, $n\geq3$, 
depends only on the associated {\em spherical} Coxeter complex.
We show moreover that it coincides with the space of 
$\De_{euc}$-valued weights 
of semistable weighted configurations on the Tits boundary $\tits X$.

The side lengths of polygons in symmetric spaces of noncompact type 
are studied in the related paper \cite{ccm}. 
Applications of the geometric results in both papers 
to algebraic group theory are given in \cite{alg}.
\end{abstract}

\tableofcontents

\section{Introduction}
For a noncompact symmetric space of rank one,
such as hyperbolic plane, 
the only isometry invariant of a geodesic segment is its metric length. 
In a symmetric space of noncompact type and arbitrary rank 
the equivalence classes of oriented segments 
modulo the identity component of the isometry group 
are parameterized by the Euclidean Weyl chamber $\Delta_{euc}$. 
We call the vector $\si(\ga)\in\De_{euc}$ corresponding to an oriented segment $\ga$ its
{\em $\De$-length}. 
The same notion of $\De$-length can be defined in a Euclidean building. 
Note that the directional part of the $\De$-length of a segment 
depends on its orientation.
This is 
because the antipodal involution of the spherical Coxeter complex 
induces an in general non-trivial involutive self-isometry of the 
spherical Weyl chamber $\De_{sph}$.
($\De_{euc}$ is the complete Euclidean cone over $\De_{sph}$.) 

For a Euclidean building or a symmetric space of noncompact type $X$
we denote by ${\cal P}_n(X)\subset \De_{euc}^n$
the set of $\De$-side lengths
which occur for oriented $n$-gons in $X$.

\begin{mthm}
\label{mthm:sidelengthsdepoosphcox}
For a thick Euclidean building $X$ 
the set ${\cal P}_n(X)$ 
of possible $\De$-side lengths of oriented $n$-gons, $n\geq3$,  
depends only on the spherical Weyl chamber $\De_{sph}$ 
associated to $X$. 
\end{mthm}

In other words, 
for any two thick Euclidean buildings 
an isomorphism of {\em spherical} Coxeter complexes,
respectively, an isometry of spherical Weyl chambers 
induces an isometry of $\De$-side length spaces. 
In particular, 
automorphisms of Coxeter complexes induce self-isometries. 

Our proof of the main theorem 
uses a relation between polygons in Euclidean buildings 
and weighted configurations on their spherical Tits buildings at infinity
via a Gauss map type construction,
see section \ref{sec:gauss}. 
A {\em weighted configuration} 
on a spherical building $B$ 
is a map 
$\psi:(\Z/n\Z,\nu)\ra B$ 
from a finite measure space. 
By composing $\psi$ with the natural projection $B\ra\De_{sph}$
onto the associated spherical Weyl chamber 
one obtains a map $(\Z/n\Z,\nu)\ra\De_{sph}$.
We call the corresponding point in $\De_{euc}^n$ 
the $\De$-{\em weights} of the configuration $\psi$.  

In order to characterize the weighted configurations on $\tits X$
which arise as Gauss maps of polygons in $X$, 
we introduce in section \ref{semi} 
a notion of {\em (semi)stability} for weighted configurations 
on abstract spherical buildings,
see also \cite[sec.\ 3.6]{ccm}.
It is motivated by Mumford stability in geometric invariant theory 
as explained in \cite[sec.\ 4]{ccm}. 
If $\psi$ is a weighted configuration on $\tits X$
then one can associate to it 
a natural convex function on $X$,
the weighted Busemann function $b_{\psi}$ 
(well-defined up to an additive constant),
and (semi)stability of $\psi$ 
amounts to certain asymptotic properties of $b_{\psi}$. 

\begin{thm}
\label{thm:samesidelengthsandweights}
Let $X$ be a Euclidean building. 
Then for $h\in\De_{euc}^n$
there exists an oriented $n$-gon in $X$ with $\De$-side lengths $h$
if and only if
there exists a semistable weighted configuration on $\tits X$ 
with $\De$-weights $h$.
\end{thm}

Balser \cite{Balser} proves the sharper result 
that the weighted configurations on $\tits X$
which arise as Gauss maps of polygons in $X$ 
are precisely the semistable ones. 

Note that every spherical building 
is the Tits boundary of a Euclidean building,
for instance, of the complete Euclidean cone over itself. 
As a step in our proof of the above results we obtain:

\begin{thm}
\label{mthm:weightsdepoosphcox}
For a thick spherical building $B$
the set of possible $\De$-weights
which occur for semistable weighted configurations
only depends on the associated spherical Weyl chamber $\De_{sph}$. 
\end{thm}

In our (logically independent) paper \cite{ccm}
we investigate the $\De$-side lengths of polygons 
in symmetric spaces of noncompact type.
We show there 
that Theorem \ref{thm:samesidelengthsandweights}
holds also in that case. 
Theorem \ref{mthm:weightsdepoosphcox} then implies 
that Theorem \ref{mthm:sidelengthsdepoosphcox} holds 
for symmetric spaces of noncompact type, too.

As a consequence of the results above and in \cite{ccm} 
it makes sense to denote by 
${\cal P}_n(\De_{sph})\subset\De_{euc}^n$ 
the space of $\De$-side lengths of oriented $n$-gons 
in {\em thick Euclidean buildings or noncompact symmetric spaces}
with spherical Weyl chamber isometric to $\De_{sph}$. 
It coincides with 
the space of $\De$-weights 
of semistable configurations on thick spherical buildings 
with this Weyl chamber. 

In most cases
the spaces ${\cal P}_n(\De_{sph})$ are known to be 
{\em finite-sided convex polyhedral cones},
namely for spherical Coxeter complexes 
which occur for a symmetric space of noncompact type.
As shown in \cite{OsheaSjamaar} and \cite{ccm} 
with rather different methods, 
${\cal P}_n(\De_{sph})$ can then be described 
as the solution set to a finite system of homogeneous linear inequalities.
The system can be given explicitely 
in terms of the Schubert calculus on Grassmann manifolds 
associated to the symmetric spaces.
The case of spherical Coxeter complexes, 
which occur for thick spherical buildings 
but not for symmetric spaces of noncompact type,
is not covered. 
(An example for such an exceptional spherical Weyl group
is the dihedral group with 16 elements $D_8$.)

\medskip
So far we discussed $\De$-side lengths.
In the case of Euclidean buildings 
there is a finer invariant for oriented geodesic segments 
taking values in $E\times E/W_{aff}$
where $(E,W_{aff})$ denotes the Euclidean Coxeter complex attached to $X$.
We call it the {\em refined} length.
Unlike the $\De$-length it also keeps track of the location of the endpoints.
In cases when the affine Weyl group acts transitively,
for examples 
for symmetric spaces 
or ultralimits of thick Euclidean buildings with cocompact affine Weyl group, 
cf.\ \cite{KleinerLeeb}, 
$\De$-length and refined length contain the same information.

An important step in our proof of the 
Main Theorem 
is a result concerning refined side lengths,
namely the observation 
that polygons can be transferred 
between thick Euclidean buildings with isomorphic {\em Euclidean} Coxeter complexes
while keeping their {\em refined} side lengths fixed,
compare Theorem \ref{prop:trans}: 
\begin{thm}
\label{introthm:trans}
For a thick Euclidean building $X$
the set ${\cal P}_n^{ref}(X)\subset(E\times E/W_{aff})^n$ 
of possible refined side lengths for $n$-gons in $X$
depends only on the associated Euclidean Coxeter complex $(E,W_{aff})$. 
\end{thm}
More generally,
polygons can be transferred to Euclidean buildings with larger affine Weyl groups 
while transforming their refined side lengths accordingly
(Addendum \ref{add:trans}),
for instance,
from a Euclidean building with one vertex 
to any other thick Euclidean building with isomorphic {\em spherical} Coxeter complex. 
Along the way
we prove analogous results for polygons in spherical buildings.

The further study of the refined side length spaces ${\cal P}_n^{ref}((E,W_{aff}))$
is relevant for certain applications to algebraic group theory
and will be taken up in \cite{alg}.

\section{Preliminaries}
\label{prelim}

In this section we will briefly review some basic facts
about singular spaces with upper curvature bound, in particular with nonpositive curvature,
and about Euclidean and spherical buildings. 
We will omit most of the proofs.
For more details on singular spaces
we refer to \cite[ch.\ 1-2]{Ballmann}, \cite[ch.\ 4+9]{Burago}, 
\cite[ch.\ 2]{KleinerLeeb} and \cite[ch.\ 2]{habil},
and for the theory of buildings from a geometric viewpoint,
i.e.\ within the framework of spaces with curvature bounded above,
to \cite[ch.\ 3-4]{KleinerLeeb}.

\subsection{Singular spaces with curvature bounded above}
\label{sec:sing}

A metric space $(Y,d)$ is called {\em geodesic} if
any two points $x,y\in Y$ can be connected by a distance minimizing geodesic
segment,
i.e.\ if there exists an isometric embedding
$\si:[0,l]\to Y$ such that $d(x,y)=l$, $\si(0)=x$ and $\si(l)=y$.
The image of such a map $\si$ is called
a {\em geodesic segment} connecting $x$ and $y$ and will be denoted by
$\ol{xy}$.
Note that this is an abuse of notation since, in general,
there may be more than one geodesic
segment connecting $x$ and $y$.


\medskip
{\bf Upper curvature bounds.}
Let $Y$ be a complete geodesic metric space.
We do {\em not} assume that $Y$ is locally compact.
One can define {\em curvature bounds} for such metric spaces
by comparison with model spaces of constant curvature.
For instance, one can compare the thickness of {\em geodesic triangles}.
Here, by a triangle we mean a one-dimensional object:
A triangle in $Y$ with the vertices $x,y,z$,
denoted by $\De=\Delta(x,y,z)$,
is the union of three geodesic segments
$\ol{xy}, \ol{yz}$ and $\ol{zx}$.
A {\em comparison triangle} $\tilde\De=\Delta(\tilde x,\tilde y,\tilde z)$ for $\Delta$
in the 2-dimensional model space $M^2_k$ with constant curvature $k$
is a triangle with the same side lengths. 
To every point $p$ on $\De$ corresponds a point $\tilde p$ on $\tilde\De$
dividing the corresponding side in the same ratio,
and we say that $\De$ is {\em thinner} than $\tilde\De$
if for any points $p,q$ on $\De$ 
the {\em chord comparison} inequality 
$d(p,q)\leq d(\tilde p,\tilde q)$ holds.
The space $Y$ has {\em curvature $\leq k$ (globally)} 
and is called a {\em $CAT(k)$-space}
if all geodesic triangles with diameter $<2\,diam(M^2_k)$
are thinner than their comparison triangles in $M^2_k$. 
In fact, if $k>0$, 
one relaxes the connectivity assumptions 
and only requires that any two points with distance $<diam(M^2_k)$
are connected by a geodesic segment. 

Due to Toponogov's Theorem,
a complete simply-connected manifold has curvature $\leq k$ 
in the distance comparison sense
if and only if it has sectional curvature $\leq k$.
A metric tree has curvature $-\infty$ in the sense that it has curvature $\leq k$ 
for all $k\in\R$.

\medskip
{\bf Angles and spaces of directions.}
The presence of a curvature bound 
allows to define {\em angles}
between geodesic segments $\si_1,\si_2:[0,\eps)\ra Y$ be unit speed geodesic segments
with the same initial point $\rho_1(0)=\rho_2(0)=y$.
Let $\tilde\al(t)$ be the angle
of a comparison triangle for $\De(y,\si_1(t),\si_2(t))$
in the appropriate model plane
at the vertex corresponding to $y$.
If $Y$ has an upper curvature bound
then the comparison angle $\tilde\al(t)$
is monotonically decreasing 
as $t\searrow0$.
It therefore converges 
and we define the angle $\angle_y(\si_1,\si_2)$ of the segments at $y$
as the limit.
In this way,
one obtains a pseudo-metric
on the space of segments emanating from a point $y\in Y$.
The metric space $(\Si_yY,\angle_y)$
obtained by identifying segments with angle zero and metric completion
is called the {\em space of directions} at $p$.
In the smooth case, $\Si_yY$ is the unit tangent sphere.
It turns out that in general $\Si_yY$ is a {\em CAT(1)-space}. 

If $\De(x,y,z)$ is a geodesic triangle and
$\tilde\De(\tilde x,\tilde y,\tilde z)$ is a comparison triangle,
the {\em angle comparison}
$\angle_x(y,z)\leq\angle_{\tilde x}(\tilde y,\tilde z)$ 
holds as a consequence of the definition of angles. 

\subsection{Hadamard spaces}
\label{sec:cat0}

We will be mainly interested in {\em CAT(0)-spaces}.
These are also called {\em Hadamard space} 
since they generalize Hadamard manifolds 
which are defined to be complete simply-connected Riemannian manifolds of nonpositive curvature.
For instance,
symmetric spaces of noncompact type are Hadamard manifolds
and Euclidean buildings are singular Hadamard spaces. 

A basic consequence of the CAT(0)-property
is the {\em convexity} of the distance function,
i.e.\ for any two constant speed geodesic segments
$\si_1,\si_2:[a,b]\ra X$ in a Hadamard space $X$ 
the distance $t\mapsto d(\si_1(t),\si_2(t))$ between fellow travellers is a convex function.
It follows that any two points can be connected by a unique geodesic segment.
In particular, Hadamard spaces are contractible.

\medskip
{\bf Boundary at infinity.}
A geodesic {\em ray} is an isometric embedding $\rho:[0,\infty)\ra X$.
By abusing notation, we will frequently identify geodesic rays
with their images.
We say that two rays are {\em asymptotic}
if they have bounded Hausdorff distance from each other
or, equivalently,
if the convex function $t\mapsto d(\rho_1(t),\rho_2(t))$ 
is bounded and hence nonincreasing. 
Asymptoticity is an equivalence relation,
and the set of equivalence classes of geodesic rays is called
the {\em ideal boundary} or {\em boundary at infinity} $\geo X$ of $X$.
An element $\xi\in\geo X$ is an {\em ideal} point or a point
{\em at infinity}.
A ray representing $\xi$ is said to be {\em asymptotic to $\xi$}.
We will use the notation $\ol{x\xi}$ to denote the unique geodesic ray
from $x\in X$ asymptotic to $\xi\in \geo X$.


The ideal boundary $\geo X$
carries a natural topology, called {\em cone topology},
which will however not play a big role in this paper. 
A basis for the cone topology is given by subsets of the following form:
For a ray $\rho_0:[0,\infty)\ra X$ and numbers $l,\eps>0$
consider all ideal points in $\geo X$ 
which are represented by rays $\rho:[0,\infty)\ra X$ such that
$d(\rho(t),\rho_0(t))<\eps$ for $0\leq t<l$.

More important for us will be a natural metric on $\geo X$,
the {\em Tits metric}. 
Given two ideal points $\xi_1,\xi_2\in\geo X$ 
we pick geodesic rays $\rho_1,\rho_2:[0,\infty)\ra X$ representing them 
and a point $x\in X$ 
and,
in analogy with the definition of angles above,
let $\tilde\al(t)$ be the angle 
of a comparison triangle in Euclidean plane 
at the vertex corresponding to $x$.
The (existence of the) limit $\lim_{t\to\infty}\tilde\al(t)$
depends only on $\xi_1,\xi_2$
and not on the location of the initial points $\rho_1(0),\rho_2(0)$
and the base point $x$.
That the limit exists follows from the observation 
that $\tilde\al(t)$ increases monotonically as $t\to\infty$
if $\rho_1(0)=\rho_2(0)=x$.
We define the {\em Tits distance} or {\em Tits angle}
$\tangle(\xi_1,\xi_2)$ to be this limit.
In other words, 
$2\sin\frac{\tangle(\xi_1,\xi_2)}{2}=
\lim_{t\to\infty}\frac{d(\rho_1(t),\rho_2(t))}{t}$.
The definition implies the useful inequality 
\begin{equation}
\label{angletangle}
\angle_x(\xi,\eta)\leq\tangle(\xi,\eta).
\end{equation}
The metric space $\tits X=(\geo X,\tangle)$
is called the {\em Tits boundary}.
As for the spaces of directions,
it turns out that the Tits boundary
is a CAT(1)-space.
Note that the Tits metric does in general {\em not} induce the cone topology.
The Tits metric is lower semicontinuous with respect to the cone
topology
and induces a topology which is (in general strictly) finer than the cone topology.

\medskip
{\bf Busemann functions.}
Busemann functions measure the relative distance from points at infinity.
They are constructed as follows. 
For an ideal point $\xi\in\geo X$
and a ray $\rho:[0,\8)\to X$ asymptotic to it
we define the {\em Busemann function} $b_{\xi}$
as the pointwise monotone limit
\[ b_{\xi}(x):=\lim_{t\to\8} (d(x,\rho(t))-t) \]
of normalized distance functions. 
One checks that, up to an additive constant,
$b_{\xi}$ does not depend on the chosen ray $\rho$.
As a limit of distance functions 
$b_{\xi}$ is convex and $1$-Lipschitz continuous. 
The level and sublevel sets of Busemann functions 
are called horospheres, respectively, horoballs. 

Convex functions have directional derivatives. 
For Busemann functions they are given by the formula
\begin{equation}
\label{busederiv}
\frac{d}{dt^+}(b_{\xi}\circ\si)(t)=-\cos\angle_{\si(t)}(\si'(t),\xi)
\end{equation}
where $\si:I\ra X$ is a unit speed geodesic segment 
and the angle on the right-hand side is taken 
between the positive direction $\si'(t)\in\Si_{\si(t)}X$ of the segment $\si$ at $\si(t)$ 
and the ray emanating from $\si(t)$ asymptotic to $\xi$.

Note that along a ray $\rho$ asymptotic to $\xi$ 
the Busemann function $b_{\xi}$ is affine linear,
i.e.\ 
$b_{\xi}(\rho(t)) = -t + const$.
As convex Lipschitz functions
Busemann functions are asymptotically linear along any ray $\rho$
and we define the {\em asymptotic slope} of $b_{\xi}$ 
at an ideal point $\eta\in\geo X$ by 
\[ slope_{\xi}(\eta)=\lim_{t\to\infty}\frac{b_{\xi}(\rho(t))}{t} \]
for a ray $\rho$ asymptotic to $\eta$.
Since $\angle_{\rho(t)}(\xi,\eta)\nearrow\tangle\xi,\eta)$ as $t\to\infty$
one obtains 
\[ slope_{\xi}(\eta)=-\cos\tangle(\xi,\eta) .\]
{\bf Cones.}
Given a metric space with diameter $\leq\pi$ 
one constructs the {\em complete Euclidean cone}
$Cone(B)$ over $B$ by mimicking the construction
which produces Euclidean 3-space from the 2-dimensional unit sphere.
The underlying set is $B\times[0,\infty)/\sim$
where $\sim$ collapses $B\times\{0\}$ to a point called the {\em tip}.
For $v_1,v_2\in B$ and $t_1,t_2\geq0$ 
we consider rays $\rho_i:[0,\infty)\to\R^2$ in Euclidean plane
with the same initial point $o$ and angle 
$\angle_o(\rho_1,\rho_2)=d_B(v_1,v_2)$.
We then define the distance of points in $Cone(B)$ 
represented by $(v_1,t_1)$ and $(v_2,t_2)$ as
$d_{\R^2}((\rho_1(t_1),\rho_2(t_2))$.

The space $Cone(B)$ is CAT(0)
if and only if $B$ is CAT(1).
In this case there is a natural isometry $B\cong\tits Cone(B)$.

\subsection{Coxeter complexes}
\label{sec:cox}

{\bf Spherical Coxeter complexes.}
A {\em spherical Coxeter complex} $(S,W_{sph})$
consists of a unit sphere $S$ and a finite subgroup $W_{sph}\subset Isom(S)$
generated by reflections.
By a reflection, we mean a reflection at a great sphere of codimension one.
$W_{sph}$ is called the {\em Weyl group} 
and the fixed point sets of the reflections in $W_{sph}$
are called {\em walls}.
The pattern of walls gives $S$ 
a natural structure of a cellular (polysimplicial) complex.
The top-dimensional cells, the {\em chambers},
are fundamental domains for the action $W_{sph}\acts S$.
They are spherical simplices if $W_{sph}$ acts without fixed point. 
If convenient,
we identify the {\em spherical model Weyl chamber} 
$\De_{sph}=S/W_{sph}$
with one of the chambers in $S$.

We will only be interested in those spherical Coxeter complexes
which are attached to semisimple complex Lie groups,
see \cite[chapter V.15]{Serre} for their classification.

An {\em embedding} $(S,W_{sph})\embed (S',W'_{sph})$
of spherical Coxeter complexes
(of equal dimensions)
consists of an isometry $\al:S\ra S'$ 
and a compatible monomorphism $\iota:W_{sph}\embed W'_{sph}$.
The isometry $\al$ maps walls to walls. 
We call the embedding of Coxeter complexes 
an {\em isomorphism}
if $\iota$ is an isomorphism.

\medskip
{\bf Euclidean Coxeter complexes.} 
A {\em Euclidean Coxeter complex} $(E,W_{aff})$
consists of a Euclidean space $E$
and a subgroup $W_{aff}\subset Isom(E)$ generated by reflections.
Again {\em reflection} means reflection at a hyperplane.
We require moreover that the induced reflection group
on the sphere $\tits E$ at infinity is finite.

One obtains an associated spherical Coxeter complex $(\tits E,W_{sph})$.
Here $W_{sph}:=rot(W_{aff})$ where $rot:Isom(E)\ra Isom(\tits E)$
is the natural homomorphism mapping
an affine transformation to its linear part.
Let $\De_{sph}$ be the spherical model Weyl chamber of $(\tits E,W_{sph})$.
We define the {\em Euclidean model Weyl chamber} 
$\De_{euc}$ of $(E,W_{aff})$
as the complete Euclidean cone over $\De_{sph}$, that is, 
$\De_{euc}=Cone(\De_{sph})$.
It is canonically identified 
with the quotient of the vector space of translations on $E$
by the natural action of $W_{sph}$ by conjugation
and one has a well-defined addition
and scalar multiplication by positive real numbers
on $\De_{euc}$. 

A  {\em wall} in the Coxeter complex $(E,W_{aff})$ is a hyperplane
fixed by a reflection in $W_{aff}$.
{\em Singular subspaces} are defined as intersections of walls,
and {\em vertices} are zero-dimensional singular subspaces.

We denote the kernel of $rot:W_{aff}\ra W_{sph}$
by $L_{trans}$ and we refer to it as the {\em translation subgroup}.
The exact sequence
$0\ra L_{trans}\ra W_{aff}\ra W_{sph}\ra 1$ splits,
i.e.\ the affine Weyl group decomposes as the semidirect product
$W_{aff} \cong W_{sph} \ltimes L_{trans}$.
The two extreme cases are that $L$ is the full group of translations on $E$,
as it happens for the Euclidean Coxeter complex attached
to a symmetric space of noncompact type,
or that $L=\{0\}$ and $W_{aff}=W_{sph}$ is finite
as in the case of Euclidean buildings with one vertex. 

A Euclidean Coxeter complex $(E,W_{aff})$ is called {\em discrete}
if $W_{aff}$ is a discrete subgroup of $Isom(E)$.
Discrete Euclidean Coxeter complexes occur as Coxeter complexes
attached to Euclidean buildings.
If moreover $W_{aff}$ acts cocompactly on $E$, 
the pattern of walls induces a 
natural structure of polysimplicial cell complex on $E$.
The top-dimensional cells, the {\em alcoves}, 
are fundamental domains for the action $W_{aff}\acts E$
and the reflections at the faces of one cell generate the group $W_{aff}$.
The alcoves are canonically isometric to the {\em model Weyl alcove}
$E/W_{aff}$. 
The Weyl alcove is different from the Euclidean Weyl chamber. 

An {\em embedding} $(E,W_{aff})\embed (E',W'_{aff})$ 
of Euclidean Coxeter complexes
(of equal dimensions) 
consists of a homothety $\al:E\ra E'$ and a compatible monomorphism
$\iota:W_{aff}\embed W'_{aff}$.
The homothety $\al$ maps walls to walls. 
We call the embedding of Coxeter complexes 
an {\em isomorphism}
if $\iota$ is an isomorphism.
A {\em dilation} of a Coxeter complex $(E,W_{aff})$ 
is a self-embedding
such that the homothety $\al:E\ra E$ is a dilation. 

\subsection{Buildings}
\label{sec:buildings}

{\bf Spherical buildings.}
A {\em spherical building} modelled on a spherical Coxeter complex $(S,W_{sph})$ 
is a CAT(1)-space $B$
together with a maximal atlas of charts,
i.e.\ isometric embeddings $S\embed B$.
The image of a chart is an {\em apartment} in $B$.
We require that any two points are contained in a common apartment 
and that the coordinate changes between charts
are induced by isometries in $W_{sph}$.

We will often denote the metric on a spherical building by $\tangle$
because in this paper spherical buildings usually 
arise as Tits boundaries.

The cell structure and the notions of wall, chamber etc.\ 
carry over from the Coxeter complex to the building. 
The building $B$ is called {\em thick}
if every codimenion-one face is adjacent to at least three chambers.
A non-thick building can always be equipped
with a natural structure of a thick building by reducing the Weyl group.
If $W_{sph}$ acts without fixed points 
the chambers are spherical simplices 
and the building carries a natural structure 
as a piecewise spherical simplicial complex. 
We will then refer to the cells as simplices. 

There is a canonical 1-Lipschitz continuous 
{\em accordion} map $acc:B\ra\De_{sph}$ 
folding the building onto the model Weyl chamber 
so that every chamber projects isometrically. 
$acc(\xi)$ is called the {\em type} of the point $\xi\in B$,
and a point in $B$ is called {\em regular}
if its type is an interior point of $\De_{sph}$.

The metric space underlying 
a $0$-dimensional spherical building modelled 
on the Coxeter complex $({\rm S}^0,\Z_2)$
is a discrete metric space where any two distinct points have distance $\pi$.

\medskip
{\bf Euclidean buildings.} 
A {\em Euclidean building} 
modelled on a Euclidean Coxeter complex $(E,W_{aff})$ 
is a CAT(0)-space $X$ 
together with a maximal atlas of charts $E\embed X$ 
subject to the following conditions:
The charts are isometric embeddings,
their images are called {\em apartments};
any pair of points and, more generally, any ray and any complete geodesic 
is contained in an apartment; 
the coordinate changes between charts are restrictions 
of isometries in $W_{aff}$.

A Euclidean building is called {\em thick}
if every wall is an intersection of apartments.
It is called {\em discrete}
if it is modelled on a discrete Euclidean Coxeter complex.
It carries then a natural structure 
as a polyhedral cell complex. 

As an example, 
the metric space underlying 
a $1$-dimensional Euclidean building modelled on $(\R,W_{aff})$,
where $W_{aff}\subset Isom(\R)$ is a sungroup generated by reflections at points,
is a metric tree.
In the discrete case it is a simplicial tree. 

If $X$ is a thick Euclidean building modelled on the Coxeter complex $(E,W_{aff})$ 
then its Tits boundary $\tits X$ is a thick spherical building
modelled on $(\tits E,W_{sph})$.
Also the spaces of directions $\Si_xX$ are spherical buildings
modelled on $(\tits E,W_{sph})$.
However, 
the building $\Si_xX$ is thick if and only if $x$ 
corresponds in a chart to a point in $E$ with maximal possible stabilizer $\cong W_{sph}$. 

If $B$ is a spherical building
then $Cone(B)$ carries a natural induced Euclidean building structure.

\section{Transfer of polygons between buildings}
\label{sec:trans}

\subsection{Polygons and side lengths}
\label{sec:lengths}

By an {\em $n$-gon} $z_1\dots z_n$
in a metric space $Z$ we mean a map $\Z/n\Z\ra Z$
carrying $i$ to the {\em vertex} $z_i$.

If $Z$ is a CAT(0)-space,
such as a Euclidean building or a symmetric space of noncompact type,
then any two points in $Z$ are connected by a unique geodesic segment
and the polygon can be promoted to a 1-dimensional object.
For any pair of successive vertices $x_{i-1}$ and $x_i$
one has a well-defined {\em side} $\ol{x_{i-1}x_i}$.
If $Z$ is a CAT(1)-space, for instance a spherical building,
one has well-defined sides for successive vertices of distance $<\pi$.
The cyclic ordering of the vertices 
determines a natural orientation of the sides. 

Let $(E,W_{aff})$ be a {\em Euclidean} Coxeter complex. 
To a pair of points $(p,q)$ in $E$ one can associate a vector 
in the Euclidean Weyl chamber $\De_{euc}=Cone(\De_{sph})$ as follows.
The translations on the affine space $E$ form a vector space
on which the spherical Weyl group $W_{sph}$ acts by conjugation. 
The quotient can be canonically identified with $\De_{euc}$.
Thus we can attach to $(p,q)$ the image in $\De_{euc}$ of the translation carrying $p$ to $q$. 
We call this vector $\si(p,q)$ 
the $\De$-{\em length} of the oriented geodesic segment $\ol{pq}$. 
It is invariant under isometries in $W_{aff}$ by construction. 
Note that the directional part of the $\De$-length 
depends on the orientation of the segment.
The reason is 
that the antipodal involution of the spherical Coxeter complex 
$(\tits E,W_{sph})$ 
induces an in general non-trivial involutive self-isometry of the 
spherical Weyl chamber $\De_{sph}$.

The complete invariant of a pair $(p,q)$
modulo the action of $W_{aff}$ is its image $\si_{ref}(p,q)$ 
under the natural projection to $E\times E/W_{aff}$. 
We call it the {\em refined length} of the oriented segment $\ol{pq}$. 
The $\De$-length is obtained by composing $\si_{ref}$ 
with the natural forgetful map
$E\times E/W_{aff}\to\De_{euc}$. 
The $\De$-length contains the complete information 
about the metric length and the direction of the segment modulo the spherical Weyl group,
while the refined length keeps track in addition of the location of the endpoints. 
If the affine Weyl group contains the full translation group,
as in the case of the Euclidean Coxeter complex attached to a symmetric space of noncompact type, 
then $E\times E/W_{aff}\cong\De_{euc}$ 
and $\De$-length and refined length contain the same information. 

As in the Euclidean case,
one can attach to a pair of points $(p,q)$ in a {\em spherical} Coxeter complex $(S,W)$ 
the {\em refined length} $\si_{ref}(p,q)\in S\times S/W$,
and it is invariant under the $W$-action.

These notions of length carry over to geometries modelled on Coxeter complexes.
One chooses an apartment containing a given pair of points 
and measures length inside the apartment. 
The length is well-defined because the coordinate changes 
between apartment charts are restrictions of isometries in the Weyl group. 

Hence one has the notion of $\De$-length 
in Euclidean buildings and symmetric spaces of noncompact type,
and one has the notion of refined length in Euclidean and spherical buildings. 
Note that in a symmetric space of noncompact type, 
although well-defined, the notion of refined length does not give more information than the $\De$-length
because the affine Weyl group acts transitively.

\medskip
Let $X$ be a Euclidean building or a symmetric space of noncompact type 
and let $(E,W_{aff})$ be its associated Euclidean Coxeter complex.
To a polygon $x_1\dots x_n$ in $X$ we associate its $\De$-side lengths 
$(\si(x_0,x_1),\dots,\si(x_{n-1},x_n))\in\De_{euc}^n$
and its refined side lengths 
$(\si_{ref}(x_0,x_1),\dots,\si_{ref}(x_{n-1},x_n))\in (E\times E/W_{aff})^n$. 
Analogously one can attach refined side lengths 
with values in $(S\times S/W)^n$
to $n$-gons in spherical buildings
modelled on the Coxeter complex $(S,W)$.

\begin{dfn}
We define ${\cal P}_n(X)\subset\De_{euc}^n$, 
respectively ${\cal P}_n^{ref}(X)\subset(E\times E/W_{aff})^n$
as the space of possible $\De$-side lengths, 
respectively refined side lengths, 
which occur for $n$-gons in $X$.
\end{dfn}

\subsection{The transfer argument}
\label{transsph}

This section is devoted to the proof of Theorem \ref{introthm:trans}
stated in the introduction. 
In fact we need to prove the same result for spherical buildings
since we will proceed by induction on the dimension 
and apply the induction assumption to the spaces of directions. 
Recall that the spaces of directions of Euclidean or spherical buildings 
are spherical buildings, 
cf.\ section \ref{sec:buildings}.

\begin{transferthm}
\label{prop:trans}
(i)
If $X$ and $X'$ are thick Euclidean buildings 
modelled on the same Euclidean Coxeter complex $(E,W_{aff})$
then ${\cal P}_n^{ref}(X)={\cal P}_n^{ref}(X')$.

(ii)
The analogous assertion for thick spherical buildings 
modelled on the same spherical Coxeter complex.
\end{transferthm}
In other words, 
isomorphisms of associated Coxeter complexes 
(cf.\ section \ref{sec:cox})
induce bijections of refined  side length spaces.
\proof
(ii)
We first discuss the spherical case.
Let $B$ and $B'$ be thick spherical buildings 
modelled on the same spherical Coxeter complex $(S,W)$.
Given a polygon $P$ in $B$, 
we will transfer it to a polygon $P'$ in $B'$
while preserving the refined side lengths. 
It suffices to show the following assertion:
{\em $(\ast)$ Let $\xi,\eta,\zeta\in B$ and $\xi',\zeta'\in B'$
so that the oriented segments $\ol{\xi\zeta}$ and $\ol{\xi'\zeta'}$
have equal refined lengths.
Then there exists $\eta'\in B'$ 
so that the triangles $\De(\xi,\eta,\zeta)$ and $\De(\xi',\eta',\zeta')$ 
have the same refined side lengths.}

We proceed by induction on the dimension of the buildings.
The assertion $(\ast)$ is trivial in dimension $0$.
We therefore assume that $dim(B)=dim(B)=d>0$ 
and that $(\ast)$ has been proven in dimensions $<d$.

There is a finite subdivision of the side $\ol{\xi\eta}$
by points $\xi_0=\xi, \xi_1,...,\xi_{k-1}, \xi_k=\eta$ such that
each geodesic triangle $\Delta(\zeta,\xi_{i},\xi_{i+1})$ is contained in an apartment.
Namely, choose the subdivision so that each subsegment $\ol{\xi_i\xi_{i+1}}$
is contained in a chamber $\Delta_i$ 
and note that each chamber $\Delta_i$ is contained in an apartment through $\zeta$.
(The analogous assertion for triangles in Euclidean buildings 
was proven in \cite[Corollary 4.6.8]{KleinerLeeb}.)

We need to find points $\xi'_1,\dots,\xi'_k$ in $B'$ 
such that the triangles $\De(\zeta',\xi'_i,\xi'_{i+1})$ 
have the same refined side lengths as $\De(\zeta,\xi_i,\xi_{i+1})$ for all $i$ 
and such that 
$\angle_{\xi'_i}(\xi'_{i-1},\xi'_{i+1})=\pi$. 
This will be done by a second induction on $i$.
We can choose $\xi'_1$ in an apartment containing $\ol{\zeta'\xi'}$.
Suppose that $\xi'_i$ has been found, $i\geq1$.
In order to find the direction $\ola{\xi'_i\xi'_{i+1}}$ at $\xi'_i$,
we apply the induction hypothesis 
(of the first induction on the dimension) 
to the links $\Si_{\xi_i}B$ and $\Si_{\xi'_i}B'$, 
which are thick spherical buildings of dimension $d-1$ 
modelled on the same spherical Coxeter complex,
and transfer the triangle 
$\De(\ola{\xi_i\xi_{i-1}},\ola{\xi_i\zeta},\ola{\xi_i\xi_{i+1}})$ in $\Si_{\xi_i}B$
to a triangle 
$\De(\ola{\xi'_i\xi'_{i-1}},\ola{\xi'_i\zeta'},\ola{\xi'_i\xi'_{i+1}})$ in $\Si_{\xi'_i}B'$
with the same refined side lengths.
We then choose an apartment in $B'$
which contains $\ol{\zeta'\xi'_i}$ 
and is tangent to the direction $\ola{\xi'_i\xi'_{i+1}}$.
Inside this apartment there is a unique choice for $\xi'_{i+1}$
with the desired properties.
After transfering all triangles $\De(\zeta,\xi_i,\xi_{i+1})$,
the concatenation of the segments $\ol{\xi'_i\xi'_{i+1}}$ 
forms a geodesic segment $\ol{\xi'\eta'}$
with the same refined length as $\ol{\xi\eta}$.
This concludes the proof in the spherical case.

(i)
The same argument works in the Euclidean case,
applying the result for spherical buildings of one dimension less.
\qed

\medskip
As a consequence, 
we can define {\em refined side length spaces}
${\cal P}_n^{ref}((E,W_{aff}))$
respectively 
${\cal P}_n^{ref}((S,W))$
associated to Euclidean and spherical Coxeter complexes. 
They describe the possible refined side lengths of polygons 
in thick buildings modelled on these Coxeter complexes. 
To be consistent with earlier notation
we may also write 
${\cal P}_n^{ref}(\De_{sph})$ instead of ${\cal P}_n^{ref}((S,W))$.

\medskip
Our proof of the Transfer Theorem \ref{prop:trans} 
allows more generally
to transfer polygons from buildings with smaller Weyl groups to buildings with larger Weyl groups. 
Suppose that 
\begin{equation}
\label{embcox}
(E,W_{aff})\embed (E',W'_{aff})
\end{equation}
is an embedding of Euclidean Coxeter complexes,
i.e.\ an isometry $E\ra E'$ 
inducing a monomorphism $W_{aff}\embed W'_{aff}$,  
cf.\ section \ref{sec:cox}.
Our above argument yields: 

\begin{add}
[to \ref{prop:trans}]
\label{add:trans}
(i)
The map 
\[ (E\times E/W_{aff})^n\to (E'\times E'/W'_{aff})^n \] 
induced by the embedding of Euclidean Coxeter complexes (\ref{embcox}) 
induces a map 
\[ {\cal P}_n^{ref}((E,W_{aff}))\lra{\cal P}_n^{ref}((E',W'_{aff})) \]
of refined side length spaces.

(ii)
The analogous assertion for embeddings of spherical Coxeter complexes. 
\end{add}
A variation of the transfer construction,
namely the folding of polygons into apartments,
will be discussed in \cite{alg}.

\section{Polygons and weighted configurations at infinity}
\label{sec:polygonconf}

The Transfer Theorem \ref{prop:trans} says 
that the possible refined side lengths for polygons 
in a thick Euclidean building
depend only on the associated affine Coxeter complex.
We now address our Main Theorem \ref{mthm:sidelengthsdepoosphcox} 
and show that the {\em unrefined} $\De$-side lengths
depend only on the {\em spherical} Coxeter complex.
That is,
we relate the $\De$-side lengths 
of polygons in Euclidean buildings with the same spherical Weyl group 
but whose affine Weyl groups may have different translation subgroups.
Addendum \ref{add:trans} 
allows to transfer polygons from buildings with smaller affine Weyl groups 
to buildings with larger ones. 
But to go in the other direction we have to pass through configurations at infinity.

We first introduce in section \ref{semi} a notion of stability 
for weighted configurations on spherical buildings 
which is motivated by 
(and consistent with, cf.\ \cite[ch.\ 4]{ccm})  
Mumford stability in geometric invariant theory. 
In section \ref{sec:gauss} we explain
how an oriented polygon in a Euclidean building $X$ 
gives rise to a collection of Gauss maps 
which can be regarded as weighted configurations 
on the spherical Tits building $\tits X$ at infinity,
and prove the basic Lemma \ref{polygonstosemistable}
that the arising configurations are semistable. 
The converse question 
when a semistable configuration $\psi$ on the Tits boundary $\tits X$
is the Gauss map of a polygon in $X$ 
amounts to a fixed point problem for a certain weak contraction $\Phi_{\psi}:X\ra X$.
In section \ref{dynamicalsystem} we prove the existence of a fixed point
in the special case when $X$ is a Euclidean building with one vertex,
i.e.\ when it is isometric to the complete Euclidean cone over its Tits boundary.
In section \ref{sec:proofsmain}
we combine our results and prove the main theorems stated in the introduction.

\subsection{Weighted configurations on spherical buildings and stability}
\label{semi}

Let $B$ be a spherical building.
We denote the metric on $B$ by $\tangle$
because spherical buildings appear in this paper usually as Tits boundaries. 

A collection 
of points $\xi_1,...,\xi_n\in B$
and of weights $m_1,...,m_n\geq0$
determines a {\em weighted configuration}
\begin{equation*}
\psi: (\Z/n\Z, \nu)\to B
\end{equation*}
on $B$.
Here $\nu$ is the measure on $\Z/n\Z$ defined by $\nu(i)=m_i$,
and the map $\psi$ sends $i$ to $\xi_i$.
By composing $\psi$ with the natural accordion projection 
$acc:B\ra\De_{sph}$
onto the associated spherical Weyl chamber $\De_{sph}$,
compare section \ref{sec:buildings}, 
one obtains a map $(\Z/n\Z,\nu)\ra\De_{sph}$.
We call the corresponding point 
$h(\psi)=(h_1,\dots,h_n)$
in $\De_{euc}^n$ 
the $\De$-{\em weights} of the configuration $\psi$,
i.e.\ $h_i=m_i\cdot acc(\xi_i)$. 
Recall that $\De_{euc}$ is defined as the complete Euclidean cone over $\De_{sph}$.

The configuration $\psi$ yields, by pushing forward $\nu$, 
the measure $\mu=\sum m_i\de_{\xi_i}$ on $B$. 
We defined its {\em slope} function on $B$ by
\begin{equation}
slope_{\mu}=-\sum_{i\in\Z/n\Z} m_i\cos\tangle(\xi_i,\cdot) .
\end{equation}
\begin{dfn}
[Stability]
The measure $\mu$ on $B$
is called {\em semistable}
if $slope_{\mu}\geq0$ 
and {\em stable}
if $slope_{\mu}>0$ everywhere on $B$. 
The weighted configuration $\psi$ is called (semi)stable
if the associated measure has this property.
\end{dfn}

\begin{ex}
Let $B$ be a spherical building of dimension $0$.
Then a measure $\mu$ on $B$ is stable
iff all atoms have mass $<\half|\mu|$,
and semistable iff all atoms have mass $\leq\half|\mu|$
(and nice semistable iff it it is either stable or consists of two atoms of the same mass).
\end{ex}

The terminology {\em slope} becomes clear when one considers 
spherical buildings as Tits boundaries. 
If $X$ is a Euclidean building 
(or a symmetric space of noncompact type)
then we can associate with a measure 
$\mu=\sum m_i\de_{\xi_i}$
on $\tits X$ 
its {\em weighted Busemann function} 
\begin{equation}
b_{\mu}:=\sum_{i\in\Z/n\Z} m_i b_{\xi_i},
\end{equation}
on $X$, 
cf.\ the definition of Busemann functions in section \ref{sec:cat0}
and the discussion of their asymptotics. 
The function $b_{\mu}$ 
is well defined up to an additive constant and {\em convex}. 
For any ideal point $\eta\in\tits X$ 
and any unit speed geodesic ray $\rho:[0,\infty)\to X$ asymptotic to $\eta$ 
holds 
\begin{equation}
slope_{\mu}(\eta)= \lim_{t\to \infty} \frac{b_\mu(\rho(t))}{t} ,
\end{equation}
i.e.\ $slope_{\mu}(\eta)$ computes the {\em asymptotic slope} of $b_\mu$ in the
direction $\eta$.

\subsection{From polygons to configurations: Gauss maps}
\label{sec:gauss}

Let $X$ be a Euclidean building or a symmetric space of noncompact type.
We now relate polygons in $X$
and weighted configurations on the spherical Tits building $\tits X$ at infinity. 

Consider a polygon $P=x_1x_2\dots x_n$ in $X$,
i.e. a map $\Z/n\Z\ra X$.
The distances $m_i=d(x_{i-1},x_i)$
determine a measure $\nu$ on $\Z/n\Z$ by putting $\nu(i)=m_i$.
The polygon $P$ gives rise to a collection $Gauss(P)$
of {\em Gauss maps}
\begin{equation}
\label{mapfromindex}
\psi:\Z/n\Z \lra \tits X
\end{equation}
by assigning to $i$ an ideal point $\xi_i\in\tits X$
so that the ray $\ol{x_{i-1}\xi_i}$ passes through $x_i$.
This construction, in the case of hyperbolic plane, 
already appears in the letter of Gauss to W.\ Bolyai \cite{Gauss}.
Taking into account the measure $\nu$,
we view the maps $\psi:(\Z/n\Z,\nu)\ra\tits X$
as {\em weighted configurations} on $\tits X$.
Their $\De$-weights equal the $\De$-side lengths of the polygon $P$. 

Note that if $X$ is a Riemannian symmetric space
and the $m_i$ are non-zero,
there is a unique Gauss map for $P$ 
because geodesic segments are uniquely extendable to complete geodesics. 
On the other hand,
if $X$ is a Euclidean building then, due to the branching of geodesics,
there are in general several Gauss maps.
However, the corresponding weighted configurations have the same $\De$-weights.


The following observation is basic for us 
and explains why the notion of semistability is useful in studying polygons.

\begin{lem}
[Semistability of Gauss maps]
\label{polygonstosemistable}
The pushed forward measures $\mu=\psi_{\ast}\nu$ are semistable.
\end{lem}
\proof
Let $\eta\in\tits X$
and let $\ga_i:[0,m_i]\ra X$ be a unit speed parametrization 
of the geodesic segment $\ol{x_{i-1}x_i}$.
Then the Busemann function $b_{\eta}$ is one-sided differentiable along $\ga_i$
with derivative
\[ 
\frac{d}{dt^+}
(b_{\eta}\circ\ga_i)(t)
= -\cos\angle_{\ga_i(t)}(\xi_i,\eta)
\leq -\cos\tangle(\xi_i,\eta) ,
\]
cf.\ the formula (\ref{busederiv}) in section \ref{sec:cat0}
for the directional derivatives of Busemann functions.
Integrating along $\ga_i$ we obtain
\begin{equation}
\label{increm}
b_{\eta}(x_i)-b_{\eta}(x_{i-1}) \leq
-m_i\cdot\cos\tangle(\xi_i,\eta) 
\end{equation}
and summation over all sides yields
\[ 0 \leq -\sum_{i\in\Z/n\Z} m_i\cdot\cos\tangle(\xi_i,\eta) = slope_{\mu}(\eta) \]
confirming the semistability.
\qed

\begin{rem}
If $X$ is a Riemannian symmetric space 
one can prove the sharper result 
that the weighted configurations on $\tits X$ 
arising as Gauss maps of closed polygons in $X$
are {\em nice} semistable,
see \cite[Lemma 5.5]{ccm}.
This refinement of the notion of semistability amounts to saying 
that the associated measures $\mu$ are semistable 
and $\{slope_{\mu}=0\}$, if non-empty, 
is a subbuilding, in fact, the Tits boundary of a totally-geodesic 
subspace of $X$,
compare \cite[Definition 3.12]{ccm} and the discussion there. 
\end{rem}

\subsection{From configurations to polygons: fixed points for weak contractions}
\label{dynamicalsystem}

Let $X$ be a Euclidean building or a symmetric space of noncompact type. 
We are now interested in finding polygons
with prescribed Gauss map.
Such polygons will correspond to the fixed points
of a certain weakly contracting self map of $X$. 

For $\xi\in \tits X$ and $t\geq0$,
we define the map 
$\phi_{\xi,t}: X\to X$
by sending $x$ to the point at distance $t$ from $x$
on the geodesic ray $\ol{x\xi}$.
Since $X$ is nonpositively curved,
the function 
$\delta: t\mapsto d(\phi_{\xi,t}(x), \phi_{\xi,t}(y))$
is convex.
It is also bounded because the rays $\ol{x\xi}$ and $\ol{y\xi}$ are asymptotic, 
and hence it is monotonically non-increasing in $t$.
This means that the maps $\phi_{\xi,t}$ are weakly contracting, 
i.e.\ they have Lipschitz constant $1$.
For a weighted configuration
$\psi:(\Z/n\Z,\nu)\ra \tits X$
we define the weak contraction 
\begin{equation}
\Phi=\Phi_{\psi}:X\lra X
\end{equation}
as the composition $\phi_{\xi_n,m_n}\circ\dots\circ\phi_{\xi_1,m_1}$.
The fixed points of $\Phi$
are the $n$-th vertices of closed polygons $P=x_1\ldots x_n$
with $\psi\in Gauss(P)$. 

Regarding the existence of fixed points for $\Phi$, 
we will only need the special case of buildings with one vertex,
that is, of complete Euclidean cones over spherical buildings. 

\begin{prop}
\label{prop:stabletofixed}
Suppose that $X$ is a Euclidean building with one vertex 
and that $\psi$ is a semistable weighted
configuration on $\tits X$.
Then the weak contraction $\Phi_\psi:X\to X$ has a fixed point.
\end{prop}

The following auxiliary result may be of independent interest.
It extends Cartan's fixed point theorem 
for isometric actions on nonpositively curved spaces with bounded orbits. 
Note that we do not need to assume local compactness.

\begin{lem}
\label{bounded}
Let $Y$ be a Hadamard space and $\Phi: Y\to Y$ a 1-Lipschitz self map.
If the forward orbits $(\Phi^n y)_{n\ge 0}$ are bounded
then $\Phi$ has a fixed point in $Y$.
\end{lem}

\proof
Consider an orbit $y_n=\Phi^n y_0$ of a point $y_0\in Y$ and
define the distance from its ``tail'' by
\[
r(y) := \limsup_{n\ra\infty} d(y_n,y).
\]
Note that $r$ inherits from the distance function the convexity
and the 1-Lipschitz continuity.
The assumption that $\Phi$ is 1-Lipschitz implies
\[
r(\Phi y)= \limsup_{n\ra\infty} d(y_n,\Phi y)
=\limsup_{n\ra\infty} d(\Phi y_{n-1},\Phi y)
\leq \limsup_{n\ra\infty} d(y_{n-1},y) =r(y),
\]
that is,
\begin{equation}
\label{phidecrr}
r\circ\Phi \leq r .
\end{equation}
It suffices to show that $r$ has a unique minimum
since this would then be a fixed point of $\Phi$.
We denote
\[ \rho:=\inf_Y r. \]
For $\eps>0$,
let $y,y'$ be points with $r(y)=r(y')<\rho+\eps$.
Then there exists $n_0$ such that for $n\geq n_0$ we have
\begin{equation*} 
\label{verticesclose}
d(y_n,y),d(y_n,y') < \rho+\eps. 
\end{equation*}
On the other hand, let $m$ be the midpoint of $\ol{yy'}$.
Since $r(m)\geq\rho$ 
we have
\begin{equation*}
\label{midpointfar}
d(y_n,m)>\rho-\eps 
\end{equation*}
for infinitely many $n$.
For these $n$ 
we apply triangle comparison to the triangle $\De(y,y',y_n)$ 
with Euclidean plane as model space,
cf.\ section \ref{sec:sing}.
For the comparison triangle $\De(\tilde y,\tilde y',\tilde y_n)$
in Euclidean plane holds the parallelogram identity: 
\[ 
d(\tilde y,\tilde y')^2+  4\, d(\tilde y_n,\tilde m)^2
= 2\bigl(d(\tilde y_n,\tilde y)^2 +
d(\tilde y_n,\tilde y')^2 \bigr)
\]
By chord comparison we have 
$d(y_n,m)\leq d(\tilde y_n,\tilde m)$ 
and obtain the inequality 
\[
d(y,y')^2+  4\, {\underbrace{d(y_n,m)}_{>\rho-\eps}}^2
\leq 2\bigl({\underbrace{d(y_n,y)}_{<\rho+\eps}}^2 +
{\underbrace{d(y_n,y')}_{<\rho+\eps}}^2 \bigr)
\]
and
\[ d(y,y')^2 < 16\rho\eps+8\eps^2 .\]
It follows that any sequence $(z_k)$ in $Y$ with $r(z_k)\searrow\rho$
must be a Cauchy sequence.
The completeness of $Y$ implies that $r$ has a minimum,
and the minimum must be unique. 
\qed

\medskip
\no
{\em Proof of Proposition \ref{prop:stabletofixed}:}
The building $X$ is isometric to the complete Euclidean cone $Cone(\tits X)$
over its Tits boundary.  
We denote the unique vertex of $X$ by $o$.

Due to the conicality of $X$ 
the contraction $\Phi$ has a fairly simple geometry. 
Let $\si$ be a simplex in $\tits X$
and let $V$ be the corresponding face of $X$,
i.e.\ the Euclidean sector with tip $o$ and ideal boundary $\si$. 
Let $\xi\in\tits X$ and $t\geq0$.
For any face $W\supseteq V$ of $X$ exists a (maximal) flat $F\subset X$
with $W\subset F$ and $\xi\in\tits F$.
The map $\phi_{\xi,t}$ restricts on $F$ to a translation 
and we have 
\begin{equation} 
\label{busincronstar}
b_{\eta}(\phi_{\xi,t}x)-b_{\eta}(x)=-t\cdot\cos\tangle(\xi,\eta) 
\end{equation}
for $\eta\in\tits F$ and $x\in F$.
Since we may vary $W$,
the equation (\ref{busincronstar}) holds for all $\eta\in\si$ and $x\in star(V)$.
We define 
$\check star(V):=\{x\in star(V):B_{|\mu|}(x)\subset star(V)\}$ 
where $|\mu|$ denotes the total mass of $|\mu|$,
i.e.\ $\check star(V)$ consists of those points in $star(V)$ 
which have at least distance $|\mu|$ from its boundary.
Since $\Phi$ has displacement $\leq|\mu|$,
we have $\Phi(\check star(V))\subset star(V)$ and 
\begin{equation}
\label{eq:busemannincrement}
b_{\eta}(\Phi x)-b_{\eta}(x)=slope_{\mu}(\eta)
\end{equation}
for $\eta\in\si$ and $x\in\check star(V)$.

We may use the Busemann functions to measure the distance from the vertex $o$.
With the normalization $b_{\eta}(o)=0$ we have
$d(o,\cdot)=\max_{\eta\in\tits X}(-b_{\eta})$. 
For technical reasons we discretize as follows.
We fix a finite subset $F\subset\De_{sph}$ in the spherical model Weyl chamber 
which has the property:
If $\eta\in\De_{sph}$ and $\zeta\in F$ such that 
$\tangle(\eta,\zeta)\leq 2\tangle(\eta,\zeta')$ for all $\zeta'\in F$
then $\eta$ has distance $>\eps$ from all faces of $\De_{sph}$
which do not contain $\zeta$ in their closure.
Let $A\subset\tits X$ be the discrete subset consisting of all points with types in $F$,
that is the inverse image of $F$ under the canonical projection $\tits X\to\De_{sph}$.
By construction $A$ satisfies:
If $\eta\in\tits X$ and $\zeta\in A$ 
such that $\tangle(\eta,\zeta)\leq 2\tangle(\eta,\zeta')$ for all $\zeta'\in A$
then $\eta$ has distance $>\eps$ from all faces of $\tits X$ 
which do not contain $\zeta$ in their closure,
i.e.\ $B_{\eps}(\eta)\subset star(\si_{\zeta})$ 
where $\si_{\zeta}$ denotes the simplex of $\tits X$ containing $\zeta$ as an interior point. 
As an approximation to $d(o,\cdot)$ we use the function
\[ f:=\max_{\zeta\in A}(-b_{\zeta}) \]
This function has bounded sublevel sets, 
and according to Lemma \ref{bounded}
we are done once we can show that $\Phi$ preserves some non-empty sublevel set. 

Let $r>0$ and $x\in X$ with $d(o,x)>r$.
For any $\zeta\in A$,
we wish to show that $-b_{\zeta}(\Phi x)\leq f(x)$. 
Suppose first that $\angle_o(x,\zeta)\leq 2\angle_o(x,\zeta')$ for all $\zeta'\in A$. 
Then, if $r$ has been chosen sufficiently large,
we have $x\in\check star(V_{\zeta})$
where $V_{\zeta}$ denotes the face of $X$ with ideal boundary $\si_{\zeta}$.
Applying (\ref{eq:busemannincrement})
and using that $\mu$ is semistable
we obtain 
$-b_{\zeta}(\Phi x)\leq -b_{\zeta}(x)\leq f(x)$
in this case.
On the other hand, 
if $\angle_o(x,\zeta)> 2\angle_o(x,\zeta')$ with $\zeta'\in A$ 
and $\angle_o(x,\zeta')=\min_{\zeta''\in A}\angle_o(x,\zeta'')$,
then 
$f(x)=-b_{\zeta'}(x)>-b_{\zeta}(x)+|\mu|$, 
again if $r$ has been chosen large enough.
So
$-b_{\zeta}(\Phi x)\leq -b_{\zeta}(x)+|\mu|\leq f(x)$
also in this case.
We conclude that $f(\Phi x)\leq f(x)$ if $d(o,x)$ is sufficiently large.
Since $\Phi$ is 1-Lipschitz
it follows that it preserves $\{f\leq R\}$ for large enough $R>0$.
\qed

\begin{rem}
\label{rem:fixedpoint}
One can show that, 
if $X$ is a Euclidean building and $\psi$ semistable 
or if $X$ is a symmetric space of noncompact type and $\psi$ nice semistable, 
then $\Phi_{\psi}$ has a fixed point.
The case of Euclidean buildings which are not necessarily locally compact
is due to Andreas Balser \cite{Balser}.
The case of symmetric spaces is proven in \cite{ccm}. 
\end{rem}

\subsection{Proofs of the main results}
\label{sec:proofsmain}

\no 
{\em Proof of  Theorem \ref{thm:samesidelengthsandweights}:}
That the existence of polygons implies the existence of configurations
follows from the semistability of Gauss maps,
cf.\ Lemma \ref{polygonstosemistable}.

The converse has been proven in Proposition \ref{prop:stabletofixed}
for Euclidean buildings with one vertex.
We deduce it for an arbitrary Euclidean building $X$ 
by using the Addendum \ref{add:trans} to the Transfer Theorem \ref{prop:trans}. 
Namely, suppose that $h$ is the $\De$-weight of a semistable configuration on $\tits X$.
Then by Proposition \ref{prop:stabletofixed}
there exist closed polygons with $\De$-side lengths $h$ 
in the Euclidean building $Cone(\tits X)$. 
Now the Euclidean Coxeter complex of $Cone(\tits X)$ embeds into the Euclidean Coxeter complex of $X$;
namely if $(E,W_{aff})$ denotes the Euclidean Coxeter complex attached to $X$
then the Euclidean Coxeter complex attached to $Cone(\tits X)$ 
is isomorphic to $(E,W_{sph})$. 
Addendum \ref{add:trans} 
therefore implies the existence of polygons with $\De$-side lengths $h$ in $X$. 
\qed

Notice that our argument does not produce a polygon having the given configuration 
as a Gauss map,
compare remark \ref{rem:fixedpoint}.

\medskip
\no
{\em Proof of   Theorem \ref{mthm:weightsdepoosphcox}:}
Consider two thick spherical buildings $B$ and $B'$ 
modelled on the same spherical Coxeter complex. 
The thick Euclidean buildings $Cone(B)$ and $Cone(B')$ 
are then modelled on the same Euclidean Coxeter complex.
The Transfer Theorem \ref{prop:trans} implies 
that in both spaces the same refined side lengths occur for closed polygons,
in particular also the same $\De$-side lengths. 
It follows from Theorem \ref{thm:samesidelengthsandweights}
that the same $\De$-weights occur for semistable weighted configurations 
on $B$ and $B'$. 
\qed

\medskip
\no
{\em Proof of Main Theorem \ref{mthm:sidelengthsdepoosphcox}:}
Consider two thick Euclidean buildings $X$ and $X'$
and suppose that their spherical Tits buildings $\tits X$ and $\tits X'$ at infinity 
are modelled on the same spherical Coxeter complex.
The spherical buildings are thick as well,
and the claim therefore follows from 
Theorems \ref{thm:samesidelengthsandweights} 
and \ref{mthm:weightsdepoosphcox}.
\qed

\medskip\no
Michael Kapovich,
Department of Mathematics,
University of California,
Davis, CA 95616, USA,
kapovich@math.ucdavis.edu

\no
Bernhard Leeb,
Mathematisches Institut,
Universit\"at M\"unchen,
Theresienstrasse 39,
D-80333 M\"unchen, Germany,
b.l@lmu.de

\no
John Millson,
Department of Mathematics,
University of Maryland,
College Park, MD 20742, USA,
jjm@math.umd.edu

\end{document}